# Algorithmic Aspects of Homeomorphism Problems.
## Alexander Nabutovsky and Shmuel Weinberger

In this paper, we will describe some results regarding the algorithmic nature of homeomorphism problems for manifolds. Very well known is the result of Novikov that for no manifold of dimension at least five is the recognition problem solvable. This is a consequence of the undecidabilty of recognition of triviality for finite presentations of groups (see section 10 of [VKP] or the appendix to [N]). For instance, by the Poincare conjecture, a PL manifold in these dimensions is the sphere iff it's homology is that of a sphere and its fundamental group is trivial. Kervaire and Milnor's [KM] work on smooth structures shows that one can also produce an algorithm to recognize the sphere among simply connected manifolds.

Indeed, there has been a certain literature on the recognition problem for restricted classes of manifolds: see e.g. [B, W1,2, Z]. Barden had speculated that modulo the problem of ascertaining whether or not fundamental groups were isomorphic, the homeomorphism problem could be decidable. However that is too optimistic as the following simple construction shows: Let $\pi$ be a group such that (1) $K(\pi,1)$ has a finite 3-skeleton and (2) $\pi$ is torsion free[0] (e.g. if the $K(\pi,1)$ is finite dimensional) and has an unsolvable word problem. (These are produced by the usual proofs of the unsolvaility of the word problem.) Let X be a finite 3-skeleton, and let M be the boundary of a regular neighborhood in $\mathbb{R}^n$ for n>6. Then $\pi_2(M) = 0$. Consider $M \# S^2 \times S^{n-3}$. $\pi_2 \cong \mathbb{Z}[\pi]$, with generator the core $S^2$. If we now surger $(g-2e)S^2$, then $\pi_2 \cong \mathbb{Z}[\pi]/(g-2e)$; so $\pi_2 \cong 0$ if g = e and contains $\mathbb{Z}[1/2]$ (by torsion freeness) if g≠e. Thus being able to decide homotopy type is the same as being able to decide diffeomorphism for thes manifolds which boils down to solving the word problem.

---
[0]This is not essential, as only a little more work shows.

Nonetheless we would like to point out that something is true in the direction suggested by Barden, namely that by combining the works of [B][S, Wi] and [GS] (and [T]) one can prove:

Theorem 1: Every PL or smooth simply connected manifold $M^n$ of dimension n at least 5 can be recognized among simply connected manifolds. That is, there is an algorithm to decide whether or not another simply connected manifold is Top, PL or Diff isomorphic to M. Moreover, an anlogous statement is true for embeddings in codimension at least three: one can algorithmically recognize any given embedding of one simply connected manifold in another up to isomorphism of pairs, or up to isotopy, if the codimension of the embedding is not two.

The proof is by concatenation of known results. Firstly for the homeomorphism problem. [Br] enables one to give a description of the first n stages of a Postnikov system for M (n=dim M). The problem is to decide whether the analogous system for some other N is homotopy equivalent. If so, obstruction theory shows that M and N are homotopy equivalent and we are most of the way done.

Deciding this can be done by induction on n. The bottom part of the Postnikov tower is the first nonzero homotopy (=homology) group. Matching up the next stage involves knowing whether the Postnikov k invariants of M and N can be moved to one another by a self homotopy equivalence of the first stage. This group of self-homotopy equivalences is clearly an arithmetic group (it is essentially $SL_n(\mathbb{Z})$ for some n) and we are asking how whether two integer vectors of some symmetric power of the defining representation lie in the same orbit. The main result of [GS] provides an algorithm for deciding whether two elements of a lattice in a vector space acted on by a given arithmetic group lie in the same orbit. For the general inductive step, note that the autohomotopy equivalences of any finite Postnikov piece form an arithmetic group by [S, Wi] so one can compare the k-invariants for the two spaces. (Note that in terms of the complexity of this proof, to complete the induction one must go from the algorithm that shows that one Postnikov stage is homotopy equivalent to another to a specific homotopy equivalence that achieves this in the next stage of comparison of k-invariants. This, in principle, can increase the complexity of the problem enormously.)

This same arithmeticity applies to the action of Aut(M) on the structure set S(M) computed by the surgery exact sequence (see [B, W3, R]). More precisely[0], one can decide first whether the rational Pontrjagin classes can be matched up, since this is a question about an action of an arithmetic group on a vector space, and then Browder-Novikov theory asserts that what remains is a finite order problem, which clearly presents no algorithmic difficulties.

In order to prove the results about embeddings one needs to understand the relevant homotopy theory and the relevant surgery. The latter can be found in Wall's book [W2] (see also [We]) where thse issues are related to the theory of "Poincare embeddings". Oddly enough, as far as we know, the relevant arithmeticity theorems for autohomotopy-equivalences of the relevant diagrams of spaces (in the sense of, say, [DK, DS]) do not appear explicitly in the literature. However, these results are implicit in [T] which deals with arithmeticity of equivariant self-homotopy equivalences. As a warm-up, let's first prove the result for self homotopy equivalences of pairs (A,B). Consider the space with involution $X = A \times E\mathbb{Z}_2 \cup B \times c(E\mathbb{Z}_2)$, where c denotes "cone"; this space is homotopy equivalent to A and has a homotopically trivial involution whose fixed point set is B. One can easily see (by obstruction theory) that the self homotopy equivalences of (A,B) form a subgroup of finite index in the equivariant self homotopy equivalences of $X^0$, which is itself arithmetc as verified in [T].

---

[0]Actually S(M) is not an affine space in the smooth category, so [GS] does not directly apply.

From an embedding, one can build an auxiliary space as above, now with $Z_6$ action whose equivariant self homotopy equivalences are commensurable with the stratified self homotopy equivalences of the pair (M, N). (Note the equivalences here are required to be stratified, so that M-N is also mapped into M-N, which is why one needs three isotropy groups, which requires the group $Z_6$ with more structure in its lattice of subgroups.) An argument similar to the case of closed manifolds completes the proof.

Remark 1: Along the way we also verified the theorem for simply connected manifolds with simply connected boundary.

Remark 2: Although group actions were used in the proof of the relevant arithmeticity statements in the theorem, one can reverse the process and use these arithmeticity statements to get information about <u>connections</u> between group actions and embeddings; see [We2]. As a simplest instance, one can use the Sullivan-Wilkerson arithmeticity results and surgery to show that the number of conjugacy classes of free actions of any finite group on an even dimensional simply connected manifold of dimension greater than four is finite. (In odd dimensions there are infinitely many conjugacy classes; finiteness is restored by adding on values of torsion and rho invariants, which are indeed variable).

---

[0] Strictly speaking, one should choose a finite CW invariant subspace of X containing A and having the same rational homotopy type. The rational version of Wall's finiteness theory easily provides such. (See e.g. [We3].)

Remark 3: By way of contrast, one can show, using the unsolvability of Hilbert's tenth problem, that the problem of the existence of a smooth embedding of M in N is, in general, unsolvable.

We now return to nonsimply connected manifolds. Since the homotopy problems are apt to be undecidable[0], it is interesting to consider instead the issue of deciding whether (simple) homotopy equivalent manifolds are isomorphic. There are two variants of the question: One asks if one can decide whether a given homotopy equivalence is homotopic to an isomorphism. The second asks whether, with the knowledge that the manifolds are homotopy equivalent, one can decide whether they are isomorphic. The trouble is that one does not specify here what the underlying homotopy equivalence should be.

We will soon see that the answer to both questions is no. The first is easier to settle. Roughly speaking the issue is that considered by surgery theory, so one needs only to concoct a situation where the surgery exact sequence can be shown uncomputable. The second problem runs into the problem of self-homotopy equivalences, and here we know nothing like arithmeticity, so conceivably the whole obstruction to isomorphism can be absorbed into some self homotopy equivalence. We will get around this by using $H\rho$ defined in [We4] for a suitable choice of antisimple manifold.

---

[0] It is, of course, very interesting to investigate how algebraicizable the homotopy theory of nonsimply connected spaces with given fundamental group is as an algorithmic "invariant" of a group. Finite groups do not present much difficulty, but even free abelian groups would reward a careful analysis.

Theorem 2: There are manifolds which cannot be distinguished algorithmically even among manifolds assumed to be (tangentially simple) homotopy equivalent to them.

We start with the work of [BDM] which produces finitely presented groups with fairly general recursively presentable group homology. The only constructions of groups used there, starting from the trivial group, are free product with amalgamation and HNN extension.

Lemma: For such groups the assembly map $A: H_*(B\pi; L) \longrightarrow L_*(\pi)$ is an isomorphism away from the prime 2.

This is very close to the result of [C] except that he only allows amalgamation along already constructed subgroups. However, if one adds in the fact that both the domain and range of the assembly map commute with direct limits, one is allowed to have arbitrary unions of constructible groups in Cappell's sense to be constructible. This slightly larger class is clearly closed under taking subgroups as "tree theory" shows: i.e. define a group of type n to be one that acts on some (perhaps) simplicial tree with type (n-1) groups as vertex and edge stabilizers. Clearly, this property is preserved under passing to subgroups, and the direct limit argument shows that Cappell's proof of the above lemma for "constructible" defined in terms of finite trees implies the version defined here.

Now take a group $\pi$ from [BDM] whose homology, while torsion free, is not algorithmic in dimension n+5. More precisely, assume that this homology group has infinitely many generators $x_i$, i=1,2,3.... and infinitely many relations, all of the form $x_i$=0 for i∈I, where I ⊂ ℕ is a recursively enumerable nonrecursive set. Let M be an antisimple manifold with fundamental group π; concretely, choose M to be the boundary of a regular neighborhood of a 2-complex with fundamental group π embedded in $\mathbb{R}^{n+1}$. By definition, Hρ(M) = 0. We can act (in the sense of the surgery exact sequence[0]) on S(M) by an arbitrary element α of $L_{n+1}(π)$ to obtain a manifold a new manifold $M_α$ tangentially simple homotopy equivalent to M. In particular we can choose α to be the A($x_n$) for n∈ℕ defining for us a sequence $M_n$. (To give an algorithmic definition of A sufficient for our purposes, one can realize elements of homology by manifolds, by a theorem of Thom, and thus constructively by a trial and error algorithm. Then one computes the symmetric signature of that manifold by a simplicial method, and rolls it up by Ranicki's algebraic surgery to be an explicit quadratic form.)

Claim:  Hρ($M_n$) = 0 iff A(n)=0 and hence iff n∈I and iff $M_n ≅ M$.

Since the homology is torsion free, $x_n$=0 iff A($x_n$) =0 iff A(h) = 0 modulo the image of $⊕H_{n+1-4i}(π; L)$, since by assumption, h ∈ $H_{n+5}$. By definition, the latter is Hρ($M_n$). Note that although Hρ depends on a choice of identification of fundamental group with π, the <u>vanishing</u> of it is independent of this, so no automorphism, even moving $π_1$ around, helps in producing an isomorphism if Hρ($M_{A(h)}$) ≠ 0.

<div style="text-align: right">QED</div>

---

[0] which essentially is an invokation of Wall's realization theorem. This is constructively proven in [W1] from a quadratic form description of L-groups.

Since the $M_n$ are constructively built manifolds tangentially simple homotopy equivalent to M, an algorithm to decide homeomorphism within this restricted class would also give an algorithm to decide membership in I, which by hypothesis is imposible, proving the theorem.

Remark 1: In this case, since A is rationally an isomorphism, and the groups involved are torsion free, there is a unique solution to the Novikov conjecture, namely the inverse of the assembly map. Thus, we could have viewed H$\rho$ as being an element of $\oplus H_{n+4i+5}(B\pi;\mathbb{Q})$, and then the argument seems even more transparent.

Remark 2: It is possible to prove the theorem using only a special case of the definition of H$\rho$, that is actually a bit simpler. We can assume that $\pi$ above is the fundamental group of an acyclic three complex (this just requires assuming that $H_1(\pi) = H_2(\pi) = 0$) and use this complex in building an antisimple homology sphere. If M' ⟶ M is a homotopy equivalence, then [M'] = [M] in $\Omega(M)$ -- where $\Omega$ denotes smooth bordism, (since $\Omega$ is a homology theory, and $\Omega(*)$ is detected by Stiefel-Whitney and Pontrjagin numbers, all of which vanish for homology spheres either because they lie in trivial groups or using the Hirzebruch signatue theorem) and hence in $\Omega(B\pi)$. This enables one to eliminate the use of intersection homology and Witt bordism theory from from [We4], which define the appropriate context for the bordism between homotopy equivalent manifolds, in general. However, one still needs the truncation algebraic nullcobordism trick from there to get a closed algebraic Poincare complex and to see the independence of the cobordism used.

Bibliography:

[B] D.Barden, Simply connected five-manifolds, Ann of Math 82 (1965) 365-385

[BDM] G.Baumslag, E. Dyer, and C.Miller, On the integral homology of finitely presented groups, Topology 23 (1983) 27-46

[Br] E.Brown, Finite computability of Postnikov complexes, Ann of Math 65 (1957) 1-20

[C] S.Cappell, On the homotopy invariance of higher signatures, Inven Math. 33 (1976) 171-179.

[DK] W.Dwyer and D.Kan, A classification theorem for diagrams of simplicial sets, Topology 23 (1984) 139-155

[DS] G.Dula and R.Schultz, Diagram cohomology and isovariant homotopy theory, Memoirs AMS 110 (1994) no 527

[GS] F.Grunewald and D.Segal, Some general algorithms, I: Arithmetic groups, Ann. of Math. 112 (1980) 531-583

[KM] M.Kervaire and J.Milnor, Groups of homotopy spheres I, Ann of Math. 77 (1963) 504-537

[N] A. Nabutovsky, Einstein structures: existence versus uniqueness, Geom. Func. Anal. 5 (1995) 76-91

[R] A.Ranicki, <u>Algebraic L-theory and Topological manifolds,</u> London Math Society, 1993

[S] D. Sullivan, Infinitesimal Computations in Topology, Publ.Math. d'IHES 47 (1977) 269-331

[T] G.Triantifillou, An algebraic model for G-homotopy types, Algebraic homotopy and local algebra (Luminy 1982) Asterisque 113-114 (1984) 312-336

[VKF] I.Volodin, Kutznetzov, and A.Fomenko, The problem of discriminating algorithmically the standard three-dimensional sphere, Russian Math. Surveys, 29:5 (1974) 71-172

[W1] C.T.C. Wall, Classification of (n-1)-connected 2n manifolds, Ann of Math 75 (1962) 163-189

[W2] ———, Classification problems in differential

topology VI: Classification of (s-1) connected 2s+1 manifolds, Topology 6 (1967) 273-296.

[W3] ———————, <u>Surgery on Compact Manifolds</u>, Academic Press 1969

[We] S.Weinberger, <u>The topologcal classification of stratified spaces</u>, University of Chicago Press 1994.

[We2] ———————, Embeddings and group actions, in preparation.

[We3] ———————, Homologically Trivial Group Actions I: Simply connected manifolds, Amer. J of Math. 108 (1986) 1005-1021

[We4] ———————, Higher rho invariants, this journal

[Wi] C.Wilkerson, Applications of minimal simplicial groups, Topology 15 (1976) 111-130 errata 19(1980) 99

[Z] A.Zubr, Classification of simply connected six manifolds, *Topology and Geometry-Rohlin Seminar*, LNM 1346 (1988) 325-339

Department of Mathematics
University of Toronto
and
Department of Mathematics
University of Chicago